\documentclass[10pt]{amsart}
\usepackage{amsfonts,amssymb}
\usepackage[all]{xy}

\newtheorem{thm}{Theorem}[section] 
\newtheorem{lem}[thm]{Lemma}

\begin{document}

\title{Admissible constants for genus 2 curves}
\author{Robin de Jong}

\thanks{The author is supported by a VENI grant 
from the Netherlands Organisation for Scientific Research (NWO). He thanks the
Max Planck Institut f\"ur Mathematik in Bonn for its hospitality during a visit.}

\begin{abstract} S.-W. Zhang recently introduced a new adelic invariant $\varphi$
for curves of genus at least~$2$ over number fields and function fields. We
calculate this invariant when the genus is equal to~$2$. 
\end{abstract}

\maketitle
\thispagestyle{empty}

\section{Introduction}

Let $X$ be a smooth projective geometrically connected curve of genus $g \geq 2$
over a field $k$ which is either a number field or the function field of a curve
over a field. Assume that
$X$ has semistable reduction over $k$. For each place $v$ of $k$, let $Nv$ be
the usual local factor connected with the product formula for $k$. 

In a recent paper \cite{zh} S.-W. Zhang proves the following theorem:
\begin{thm} Let $(\omega,\omega)_a$ be the admissible self-intersection of the
relative dualizing sheaf of $X$. Let $\langle
\Delta_\xi,\Delta_\xi \rangle$ be the height of the canonical Gross-Schoen cycle
on $X^3$. Then the formula:
\[ (\omega,\omega)_a = \frac{2g-2}{2g+1} \left( \langle \Delta_\xi,\Delta_\xi
\rangle + \sum_v \varphi(X_v) \log Nv \right) \]
holds, where the $\varphi(X_v)$ are local invariants associated to 
$X \otimes k_v$, defined as follows:
\begin{itemize}
\item if $v$ is a non-archimedean place, then:
\[ \varphi(X_v) = -\frac{1}{4} \delta(X_v) + \frac{1}{4} \int_{R(X_v)}
g_v(x,x)((10g+2) \mu_v - \delta_{K_{X_v}}) \, , \]
where:
\begin{itemize}
\item $\delta(X_v)$ is the number of singular points on the special
fiber of $X \otimes k_v$,
\item $R(X_v)$ is the reduction graph of $X\otimes k_v$, 
\item $g_v$ is the
Green's function for the admissible metric $\mu_v$ on $R(X_v)$,   
\item $K_{X_v}$ is the
canonical divisor on $R(X_v)$. 
\end{itemize} 
In particular, $\varphi(X_v)=0$ if $X$ has good reduction at $v$;
\item if $v$ is an archimedean place, then:
\[ \varphi(X_v) = \sum_{\ell} \frac{2}{\lambda_\ell} \sum_{m,n=1}^g  \left| 
\int_{X(\bar{k}_v)} \phi_\ell \omega_m \bar{\omega}_n \right|^2 \, , \]
where $\phi_\ell$ are the normalized real eigenforms of the Arakelov 
Laplacian on
$X(\bar{k}_v)$ with eigenvalues 
$\lambda_\ell>0$, and $(\omega_1,\ldots,\omega_g)$ is an orthonormal basis for 
the hermitian inner
product $(\omega,\eta) \mapsto \frac{i}{2} \int_{X(\bar{k}_v)} \omega \, \bar{\eta}$
on the space of holomorphic differentials.
\end{itemize}
\end{thm}
Apart from giving an explicit connection between the two canonical invariants
$(\omega,\omega)_a$ and $\langle \Delta_\xi,\Delta_\xi \rangle$,
Zhang's theorem has a possible application to the effective Bogomolov
conjecture, \emph{i.e.}, the question of giving effective positive lower bounds for
$(\omega,\omega)_a$. Indeed, the height of the canonical Gross-Schoen cycle 
$\langle \Delta_\xi,\Delta_\xi \rangle$ is known to be non-negative in the case
of a function field in characteristic zero, and should be non-negative in
general by a standard conjecture of Gillet-Soul\'e 
(\emph{op. cit.}, Section~2.4). 
Further, the invariant $\varphi$ should be non-negative, and Zhang proposes, 
in the non-archimedean case, an
explicit lower bound for it which is positive in the case of non-smooth
reduction (\emph{op. cit.}, Conjecture~1.4.2). 
Note that it is clear from the definition that $\varphi$ is
non-negative in the archimedean case; in fact it is positive (\emph{op. cit.},
Remark after Proposition~2.5.3).

Besides $\varphi(X_v)$, Zhang also considers the invariant
$\lambda(X_v)$ defined by:
\[ \lambda(X_v) = \frac{g-1}{6(2g+1)} \varphi(X_v) + \frac{1}{12} (\varepsilon(X_v)
+ \delta(X_v)) \, , \]
where:
\begin{itemize}
\item if $v$ is a non-archimedean place, the invariant $\delta(X_v)$ is as
above, and:
\[ \varepsilon(X_v) = \int_{R(X_v)} g_v(x,x)((2g-2)\mu_v + \delta_{K_{X_v}} )\, ,
\]
\item if $v$ is an archimedean place, then:
\[ \delta(X_v) = \delta_F(X_v) - 4g \log(2\pi) \]
with $\delta_F(X_v)$ the Faltings delta-invariant of the compact Riemann surface
$X(\bar{k}_v)$, and $\varepsilon(X_v)=0$. 
\end{itemize}
The significance of this invariant is that if $\deg \det R\pi_* \omega$ 
denotes the
(non-normalized) geometric or Faltings height of $X$ one has a simple expression:
\[ \deg \det R\pi_* \omega = \frac{g-1}{6(2g+1)} \langle \Delta_\xi,\Delta_\xi \rangle + \sum_v
\lambda(X_v) \log Nv  \]
for $\deg \det R\pi_*\omega$, as follows from the Noether formula:
\[ 12 \deg \det R\pi_* \omega = (\omega,\omega)_a + \sum_v
(\varepsilon(X_v)+\delta(X_v)) \log Nv \, . \]

Now assume that $X$ has genus $g=2$. Our purpose is to calculate the invariants
$\varphi(X_v)$ and $\lambda(X_v)$ explicitly. For the $\lambda$-invariant we
obtain:
\begin{itemize}
\item if $v$ is non-archimedean, then:
\[ 10 \lambda(X_v) = \delta_0(X_v) + 2\delta_1(X_v) \, , \]
where $\delta_0(X_v)$ is the number of non-separating
nodes and $\delta_1(X_v) $ is the number of separating nodes 
in the special fiber of $X \otimes k_v$;
\item if $v$ is archimedean, then: 
\[ 10 \lambda(X_v) = -20 \log(2\pi) - \log \|\Delta_2\|(X_v) \, , \]
where $\|\Delta_2\|(X_v)$ is the normalized modular
discriminant of the compact Riemann surface $X(\bar{k}_v)$ (see
below). 
\end{itemize}
Thus, the $\lambda(X_v)$ are precisely the well-known local invariants
corresponding to the discriminant modular form of weight~$10$ \cite{li}
\cite{sa} \cite{ue}. In particular we have:
\[ \deg \det R\pi_* \omega = \sum_v \lambda(X_v) \log Nv \]
and we recover the fact that the height of the canonical Gross-Schoen
cycle vanishes for $X$.

\section{The non-archimedean case}

Let $k$ be a complete discretely valued field. Let $X$ be a smooth projective
geometrically connected curve of genus~$2$ over $k$. Assume that $X$ has
semistable reduction over $k$. In this section we give
the invariants $\varphi(X)$ and $\lambda(X)$ of $X$. 

The proof of our result is based on the classification of the semistable
fiber types in genus~$2$ and consists of a case-by-case analysis. The notation
we employ for the various fiber types is as in \cite{mo2}.
We remark that there are
no restrictions on the residue characteristic of $k$.
\begin{thm} The invariant $\varphi(X)$ is given by the following table,
depending on the type of the special fiber of the regular minimal model of 
$X$:
\begin{center}
\begin{tabular}{|l|cccc|}
\hline  
Type & $\delta_0$ & $\delta_1$ & $\varepsilon$ & $\varphi$ \\  [2pt] 
\hline  
$I$  & $0$ & $0$ & $0$ & $0$ \\  [4pt]
$II(a)$ & $0$ & $a$ & $a$ & $a$ \\  [4pt] 
$III(a)$ & $a$ & $0$ & $\frac{1}{6}a$ & $\frac{1}{12}a $ \\ [4pt]
$IV(a,b)$ & $b$ & $a$ & $a+\frac{1}{6}b$ & $a + \frac{1}{12}b$ \\ [4pt]
$V(a,b)$ & $a+b$ & $0$ & $\frac{1}{6}(a+b)$ & $ \frac{1}{12}(a+b)$ \\ [4pt]
$VI(a,b,c)$ & $b+c$ & $a$ & $a+ \frac{1}{6}(b+c)$ & $a + \frac{1}{12}(b+c) $ \\
[4pt]
$VII(a,b,c)$ & $a+b+c$ & $0$ & 
$\frac{1}{6}(a+b+c) + \frac{1}{6} \frac{abc}{ab+bc+ca}$ &
$\frac{1}{12}(a+b+c) - \frac{5}{12}\frac{abc}{ab+bc+ca}$ \\  [4pt]
\hline
\end{tabular}
\end{center}
\vspace{8pt}
For $\lambda(X)$ the formula:
\[ 10 \lambda(X) = \delta_0(X) + 2 \delta_1(X) \]
holds.
\end{thm}
Let us indicate how the theorem is proved. Let
$r$ be the effective resistance function on the reduction graph $R(X)$ of $X$, 
extended bilinearly to a pairing on $\mathrm{Div}(R(X))$. By Corollary~2.4 of
\cite{fa} the formula:
\[ \varphi(X) = -\frac{1}{4}(\delta_0(X) + \delta_1(X)) -\frac{3}{8}r(K,K) + 2
\varepsilon(X)  \]
holds, where $K$ is the canonical divisor on $R(X)$. 
The invariant $r(K,K)$ is calculated by viewing $R(X)$ as an electrical circuit.
The invariant $\varepsilon$  
is calculated on the basis of explicit expressions for 
the admissible measure and admissible Green's function; see \cite{mo1} and 
\cite{mo2} for such computations. The results we find are as follows:
\begin{center}
\begin{tabular}{|l|cccc|}
\hline  
\emph{Type} & $\delta_0$ & $\delta_1$ & $r(K,K)$ & $\varepsilon$ \\  [2pt] 
\hline  
$I$  & $0$ & $0$ & $0$ & $0$ \\  [4pt]
$II(a)$ & $0$ & $a$ & $2a$ & $a$ \\  [4pt] 
$III(a)$ & $a$ & $0$ & $0$ & $\frac{1}{6}a $ \\ [4pt]
$IV(a,b)$ & $b$ & $a$ & $2a$  & $a+\frac{1}{6}b$ \\ [4pt]
$V(a,b)$ & $a+b$ & $0$ & $ 0$ & $\frac{1}{6}(a+b)$ \\ [4pt]
$VI(a,b,c)$ & $b+c$ & $a$ &   $2a $ & $a+ \frac{1}{6}(b+c)$\\
[4pt]
$VII(a,b,c)$ & $a+b+c$ & $0$ & 
$2\frac{abc}{ab+bc+ca}$ & 
$\frac{1}{6}(a+b+c) + \frac{1}{6} \frac{abc}{ab+bc+ca}$\\  [4pt]
\hline
\end{tabular}
\end{center}
\vspace{8pt}
The values of $\varphi$ follow. 

The formula for $\lambda(X)$ is verified for each case separately.

\section{The archimedean case}

Let $X$ be a compact and connected Riemann surface of genus~$2$. 
In this section
we calculate the invariants $\varphi(X)$ and $\lambda(X)$ of $X$. Let
$\mathrm{Pic}(X)$ be the Picard variety of $X$, and for each integer~$d$ 
denote by $\mathrm{Pic}^d(X)$ the component of $\mathrm{Pic}(X)$ of 
degree~$d$. We have a canonical 
theta divisor $\Theta$ on $\mathrm{Pic}^1(X)$, and a standard hermitian 
metric $\| \cdot \|$ on the line 
bundle $\mathcal{O}(\Theta)$ on
$\mathrm{Pic}^1(X)$. Let $\nu$ be its curvature form. We have:
\[  \int_{\mathrm{Pic}^1(X)} \nu^2 = \Theta . \Theta = 2 \, . \]
Let $K$ be a canonical divisor on $X$, and let $\mathbf{P}$ be the set of~$10$ points $P$ 
of $\mathrm{Pic}^1(X) - \Theta$ such that $2P \equiv K$. Denote by
$\|\theta\|$ the norm of the canonical section $\theta$ of 
$\mathcal{O}(\Theta)$. We let:
\[ \|\Delta_2\|(X) = 2^{-12} \prod_{P \in \mathbf{P}} \|\theta\|^2(P) \, , \]
the normalized modular discriminant of $X$, and we let $\|H\|(X)$ be the invariant
of $X$ defined by:
\[ \log \|H\|(X) = \frac{1}{2} \int_{\mathrm{Pic}^1(X)} \log \|\theta\| \, \nu^2 
\, . \]
These two invariants were introduced in \cite{bo}.
\begin{thm} For the $\varphi$-invariant and the $\lambda$-invariant of $X$, the
formulas: 
\[
\varphi(X) = -\frac{1}{2} \log \|\Delta_2\|(X) + 10 \log \|H\|(X) \]
and
\[ 10 \lambda(X)  =  -20 \log(2\pi) - \log \|\Delta_2 \|(X) \]
hold.
\end{thm}
The key to the proof is the following lemma. Let $\Phi$ be the map:
\[ X^2 \to \mathrm{Pic}^1(X) \, , \qquad (x,y) \mapsto [2x - y] \, . \]
\begin{lem} \label{uniquelemma} 
The map $\Phi$ is finite flat of degree~$8$.
\end{lem}
\begin{proof} Let $y \mapsto y'$ be the hyperelliptic involution of $X$. We have
a commutative diagram:
\[
\xymatrix{ X^2 \ar[d]_\alpha \ar[r]^\Phi & \mathrm{Pic}^1(X)  \\
X^2 \ar[r]^{\Phi^\lor} & \mathrm{Pic}^3(X) \ar[u]_\beta
}
\]
where $\alpha$ and $\beta$ are isomorphisms, with:
\[ \begin{array}{c}
\alpha \colon X^2 \to X^2 \, , \\
(x,y) \mapsto (x,y') \, , \end{array}  \quad 
\begin{array}{c}
\Phi^\lor \colon X^2 \to \mathrm{Pic}^3(X) \, , \\
(x,y) \mapsto [2x+y] \, , \end{array}  \quad  
 \begin{array}{c}
\beta \colon \mathrm{Pic}^3(X) \to \mathrm{Pic}^1(X) \, , \\
\, [ D ] \mapsto [D - K] \, . \end{array}   
\]
It suffices to prove that $\Phi^\lor$ is finite flat of
degree~$8$. Let $p \colon X^{(3)} \to \mathrm{Pic}^3(X)$ be the natural map;
then $p$ is a $\mathbf{P}^1$-bundle over $\mathrm{Pic}^3(X)$, and $\Phi^\lor$
has a natural injective lift to $X^{(3)}$. A point $D$ on $X^{(3)}$ is in the
image of this lift if and only if $D$, when seen as an effective divisor on $X$, 
contains a point which is ramified for the morphism
$X \to \mathbf{P}^1$ determined by the fiber $|D|$
of $p$ in which $D$ lies. Since every morphism 
$X \to \mathbf{P}^1$ associated to a 
$D$ on $X^{(3)}$ is ramified, the map $\Phi^\lor$ is surjective. 
As every morphism $X \to \mathbf{P}^1$ associated to a $D$
on $X^{(3)}$ has only finitely many ramification points, the map $\Phi^\lor$ is
quasi-finite, hence finite since $\Phi^\lor$ is proper. As $X^2$ and $\mathrm{Pic}^3(X)$
are smooth and the fibers of $\Phi^\lor$ are equidimensional, 
the map $\Phi^\lor$ is flat. By 
Riemann-Hurwitz the generic $X \to \mathbf{P}^1$ associated to a 
$D$ on $X^{(3)}$ has~$8$
simple ramification points. It follows that the degree of $\Phi^\lor$ is~$8$.
\end{proof}
Let $G \colon X^2 \to \mathbf{R}$ be the Arakelov-Green's function of $X$, and
let $\Delta$ be the diagonal divisor on $X^2$. We have a canonical hermitian
metric on the line bundle $\mathcal{O}(\Delta)$ on $X^2$ by putting $\|1\|(x,y)
= G(x,y)$, where~$1$ is the canonical section of $\mathcal{O}(\Delta)$. Denote
by $h_\Delta$ the curvature form of $\mathcal{O}(\Delta)$. We have:
\[ \int_{X^2} h_\Delta^2 = \Delta . \Delta = -2 \, . \]
Restricting $\mathcal{O}(\Delta)$ to a fiber of any of the two natural 
projections of $X^2$ onto $X$
and taking the curvature form we obtain the Arakelov $(1,1)$-form $\mu$ on $X$.
We have $\int_X \mu =1$ and:
\[ \int_X \log G(x,y) \, \mu(x) = 0 \]
for each $y$ on $X$.
Let $(\omega_1,\omega_2)$ be
an orthonormal basis of $\mathrm{H}^0(X,\omega_X)$, the space of holomorphic
differentials on $X$. We can write explicitly:
\[ h_\Delta(x,y) = \mu(x) + \mu(y) - i \sum_{k=1}^2 (\omega_k(x)
\bar{\omega}_k(y) + \omega_k(y) \bar{\omega}_k(x) )   \]
and:
\[ \mu(x) = \frac{i}{4} \sum_{k=1}^2 \omega_k(x) \bar{\omega}_k(x) \, . \]
By \cite[Proposition 2.5.3]{zh} we have:
\[ \varphi(X) = \int_{X^2} \log G \, h_\Delta^2 \, . \]
We compute the integral using our results from \cite{dj1} and
\cite{dj2}. 
Let $W$ be the divisor of Weierstrass points on $X$, and let $p_1
\colon X^2 \to X$ be the projection onto the first coordinate. The divisor $W$ is
reduced effective of degree~$6$. 
According to \cite[p.~31]{fay} there exists a canonical isomorphism:
\[ \sigma \colon \Phi^* \mathcal{O}(\Theta) \xrightarrow{\cong} \mathcal{O}(2\Delta + p_1^*W)
\] 
of line bundles on $X^2$, identifying the canonical sections on both sides. 
In \cite[Proposition 2.1]{dj1} we proved that this
isomorphism has a constant norm over $X^2$. Thus, the curvature forms on both
sides are equal:
\[ \Phi^* \nu = 2 h_\Delta + 6 \mu(x)  \quad \textrm{on} \,\, X^2 \, . \]
Squaring both sides of this identity we get:
\[ h_\Delta^2 = \frac{1}{4} \Phi^*(\nu^2) - 6 h_\Delta \mu(x) \, , \]
since $\mu(x)^2=0$.
Denote by $S(X)$ the norm of $\sigma$. Then we have:
\[ 2 \log G(x,y) + \sum_w \log G(x,w) = \log \|\theta\|(2x-y) + \log S(X) \]
for generic $(x,y) \in X^2$, where $w$ runs through the Weierstrass points of $X$.
By fixing $y$ and integrating against $\mu(x)$ on $X$ we find that:
\[ \log S(X) = - \int_X \log\|\theta\|(2x-y) \, \mu(x) \, . \]
By integrating against $h_\Delta^2$ on $X^2$ we obtain:
\[ 2 \varphi(X) + \sum_w \int_{X^2} \log G(x,w) \, h_\Delta^2 = -2 \log S(X)
+ \int_{X^2} \log \|\theta\|(2x-y) \, h_\Delta^2 \, . \]
As we have:
\[ h_\Delta^2 = 2 \mu(x)\mu(y) - \sum_{k,l=1}^2 (\omega_k(x) \bar{\omega}_l(x)
\bar{\omega}_k(y)\omega_l(y) +
\bar{\omega}_k(x)\omega_l(x)\omega_k(y)\bar{\omega}_l(y)) \]
it follows that:
\[ \int_{X^2} \log G(x,w) \, h_\Delta^2 = 0 \]
for each $w$ in $W$ and hence we simply have:
\[ 2 \varphi(X) = -2 \log S(X)
+ \int_{X^2} \log \|\theta\|(2x-y) \, h_\Delta^2 \, . \]
Using our earlier expression for $h_\Delta^2$ this becomes:
\[ 2\varphi(X) = -2 \log S(X) + \int_{X^2} \log \|\theta\|(2x-y) \left(
\frac{1}{4} \Phi^*(\nu^2) - 6 h_\Delta \mu(x) \right) \, . \]
It is easily verified that $h_\Delta \mu(x) =h_\Delta \mu(y) =  \mu(x) \mu(y)$ and
hence:
\[ \int_{X^2} \log \|\theta\|(2x-y) \, h_\Delta \mu(x) =
\int_{X^2} \log \|\theta\|(2x-y) \, \mu(x) \mu(y) = -\log S(X) \, . \]
From Lemma \ref{uniquelemma} it follows that:
\[ \int_{X^2} \log \|\theta\|(2x-y) \, \Phi^*(\nu^2) = 8 \int_{\mathrm{Pic}^1(X)} 
\log \|\theta\| \, \nu^2 = 16 \log \|H\|(X) \, . \]
All in all we find:
\[ \varphi(X) = 2 \log S(X) + 2 \log \|H\|(X) \, . \]
Let $\delta_F(X)$ be the Faltings delta-invariant of $X$. According to
\cite[Corollary 1.7]{dj2} the formula:
\[ \log S(X) = -16 \log(2\pi) - \frac{5}{4} \log \|\Delta_2\|(X) - \delta_F(X)
 \]
holds, and in turn, according to \cite[Proposition 4]{bo} we have:
\[ \delta_F(X) = -16 \log(2\pi) - \log \|\Delta_2\|(X) - 4 \log \|H\|(X) \, . \]
The formula:
\[ \varphi(X) = -\frac{1}{2} \log \|\Delta_2\|(X) + 10 \log \|H\|(X) \]
follows. 

By definition we have:
\[ \lambda(X) = \frac{1}{30} \varphi(X) +\frac{1}{12} \delta_F(X) -
\frac{2}{3} \log(2\pi) \]
so we obtain:
\[ 10 \lambda(X) =  -20 \log(2\pi) - \log \|\Delta_2 \|(X) \]
by using \cite[Proposition 4]{bo} once more.

\vspace{1cm}

\noindent Address of the author: \\  \\
Robin de Jong \\
Mathematical Institute \\
University of Leiden \\
PO Box 9512 \\
2300 RA Leiden \\
The Netherlands \\
Email: \verb+rdejong@math.leidenuniv.nl+

\end{document}